\documentclass[a4paper,12pt]{amsart}
\usepackage{amsmath,amsthm,amssymb,latexsym,epic}
\usepackage{graphicx,enumerate}

\newtheorem{theorem}{Theorem}
\newtheorem{lemma}[theorem]{Lemma}
\newtheorem{remark}[theorem]{Remark}
\newtheorem{corollary}[theorem]{Corollary}
\newtheorem{proposition}[theorem]{Proposition}

\newtheorem{problem}[theorem]{Problem}
\usepackage[all]{xy}

\newcommand{\FP}{{\mathcal{FP}^+(S_n)}}
\newcommand{\fp}{{\mathcal{FP}^+}}
\begin{document}
\title{Simple modules over factorpowers}
\author{Volodymyr Mazorchuk}
\date{}

\maketitle

\begin{abstract}
In this paper we study complex representations of the factorpower
$\fp(G,M)$ of a finite group $G$ acting on a finite set $M$. This
includes the finite monoid $\FP$, which can be seen as a kind of 
a ``balanced'' generalization of the symmetric group $S_n$ inside 
the semigroup of all binary relations. We describe all irreducible representations of $\fp(G,M)$ and relate them to irreducible 
representations of certain inverse semigroups. In particular, 
irreducible representations of $\FP$ are related to irreducible 
representations of the maximal factorizable submonoid of the dual 
symmetric inverse monoid. We also show that in the latter cases
irreducible representations lead to an interesting combinatorial problem 
in the representation theory of $S_n$, which, in particular, is 
related to Foulkes' conjecture. Finally, we show that all simple
$\fp(G,M)$-modules are unitarizable and that tensor products of simple
$\fp(G,M)$-modules are completely reducible.
\end{abstract}

\section{Introduction and preliminaries}\label{s1}

Let us first fix some notation. For a semigroup $S$ we denote by
$E(S)$ the set of idempotents of $S$; and by $\mathcal{D}$,
$\mathcal{L}$, $\mathcal{R}$, $\mathcal{H}$, $\mathcal{J}$ 
the corresponding Green's relations on $S$. 

Let $S$ be a finite semigroup acting on a finite set $M$ by (everywhere
defined) transformations. Consider the power semigroup $\mathcal{P}(S)$, 
which consists of all subsets of $S$ with the natural multiplication
$A\cdot B=\{a\cdot b|a\in A, b\in B\}$ for $A,B\in \mathcal{P}(S)$. 
Define on $\mathcal{P}(S)$ the binary relation $\sim_M$ in the following
way: $A\sim_M B$ provided that $\{a(m)|a\in A\}=\{b(m)|b\in B\}$. 
It is straightforward to verify that $\sim_M$ is a congruence on 
$\mathcal{P}(S)$ and hence we can consider the corresponding quotient
$\mathcal{FP}(S,M)=\mathcal{P}(S)/_{\sim_M}$. The semigroup 
$\mathcal{FP}(S,M)$ has an isolated zero element, namely the class of
the empty set. Denote by $\mathcal{FP}^+(S,M)$ the complement to this
zero element. The semigroup $\mathcal{FP}^+(S,M)$ is called the 
{\em factorpower} of the action of $S$ on $M$. This construction was 
proposed and first studied in \cite{GM1}. In general, there might exist 
many different nonequivalent actions of $S$ on different sets, which give 
rise to different (nonisomorphic) factor powers. In \cite{D1,D2} all 
factorpowers of a finite group $G$ are classified up to isomorphism. 
In the present paper we would like to study complex representations of
these factorpowers.

There is one very special factorpower, which is studied much better than
all others and has some really exciting properties. For a positive integer
$n$ let $\FP$ denote the factorpower of the symmetric 
group $S_n$ with respect to its natural action on the set 
$\mathbf{N}=\mathbf{N}_n=\{1,2,\dots,n\}$. This semigroup was studied
in \cite{GM1,GM2,GM3,GM4,M1,M2}. From the
definition it follows that $\mathcal{FP}^+(S,M)$ is a subsemigroup  of
the semigroup $\mathfrak{B}(M)$ of all binary relations on $M$, in
particular, $\FP$ is a subsemigroup  of $\mathfrak{B}_n$, the semigroup
of all binary relations on $\mathbf{N}$. The semigroup $\FP$ has
the following properties (see \cite{GM2}): 
\begin{enumerate}[(I)]
\item\label{enu-1} Idempotents of $\FP$ are exactly the equivalence 
relations on $\mathbf{N}$, in particular, asymptotically almost all 
idempotents of $\mathfrak{B}_n$ do not belong to $\FP$ in the sense
$\frac{|E(\FP)|}{|E(\mathfrak{B}_n)|}\to 0$, $n\to \infty$.
\item\label{enu-2} The semigroup $\FP$ contains asymptotically almost 
all elements of $\mathfrak{B}_n$ in the sense that 
$\frac{|\FP|}{|\mathfrak{B}_n|}\to 1$, $n\to \infty$.
\item\label{enu-3} The semigroup $\FP$ is the maximum  subsemigroup 
of $\mathfrak{B}_n$ which contains $S_n$ and whose idempotents are 
exactly the equivalence relations.
\item\label{enu-4} The semigroup $\FP$ is a natural quotient of the 
semigroup of doubly-stochastic $n\times n$ matrices.
\item\label{enu-5} $\FP$ is a universal factorpower in the sense that
every factorpower of a finite group acting on a finite set is a 
subsemigroup of $\FP$ for some $n$.
\end{enumerate}
In \cite{GM3,M1,M2} one can find classifications of maximal 
nilpotent subsemigroups and automorphisms of $\FP$, as well as a
description of Green's relations.

Being such a natural generalization of the symmetric group one could
expect $\FP$ to have further interesting properties. The motivation for
the present paper was an attempt to understand complex representations 
of $\FP$. There is a well-developed abstract theory for the study of
irreducible complex representations of semigroups, see \cite{Mu,CP,GMS,GM5}.
Trying to apply it to the study of complex representations of $\FP$
revealed that $\FP$ and, more generally, factorpowers of finite groups
have some nice properties. For example, they all have an involution and
all regular $\mathcal{D}$-classes of such factorpowers are in fact inverse.
This observation allows for a straightforward application of the general
theory, which results in a nice description of irreducible
representations of $\fp(G,M)$, where $G$ is a  group, in terms of 
certain induced $G$-modules. This is worked out in Section~\ref{s2}. 
As a corollary we obtain that the quotient of $\mathbb{C}[\FP]$ modulo 
the Jacobson ideal is isomorphic to the semigroup algebra of the maximal 
factorizable submonoid $\mathcal{F}_n^*$ of the  dual symmetric inverse 
monoid $\mathcal{I}_n^*$ which was introduced in \cite{FL}. In particular,
irreducible representations of $\FP$ and $\mathcal{F}_n^*$ are the same
(at least as $S_n$-modules, but also as modules over the traces of the
corresponding $\mathcal{D}$-classes).

It further turns out that the representation theory of $\FP$ (and also 
that of $\mathcal{FI}_n^*$) has an interesting relation to the 
representation theory of the symmetric group $S_n$. The group $S_n$ is 
the group of units of $\FP$. In particular, every representation of 
$\FP$ becomes a representation of $S_n$ via restriction. In the study of 
representations of semigroups a special role is played by modules, 
associated with regular $\mathcal{L}$-classes. Every regular 
$\mathcal{L}$-class of $\FP$ contains a unique idempotent. This idempotent 
is an equivalence relation on $\mathbf{N}$ and the sizes of the 
blocks of this relation determine a partition of $n$, say $\lambda$. 
We will show that, after restriction to $S_n$, the representation of 
$\FP$ associated with our $\mathcal{L}$-class is isomorphic to the 
permutation module $\mathcal{M}^{\lambda}$ (see \cite[2.1]{Sa}). The module
$\mathcal{M}^{\lambda}$ is usually realized via tabloids of shape 
$\lambda$. If these tabloids have rows of the same length, permutations 
of such rows induce automorphisms of $\mathcal{M}^{\lambda}$, which give 
rise to the action of a certain group $G$ on $\mathcal{M}^{\lambda}$ 
via automorphisms. This group $G$ turns out to be the (unique) maximal 
subgroup of $\FP$, contained in our $\mathcal{L}$-class. Considering 
the isotypic components of $\mathcal{M}^{\lambda}$ with respect to the 
action of $G$ leads to a finer decomposition of $\mathcal{M}^{\lambda}$ 
as an $S_n$-module. The elements of this decomposition are simple 
$\FP$-modules and we raise the problem of determining the multiplicities 
of Specht modules in these simple $\FP$-modules. For the permutation module 
$\mathcal{M}^{\lambda}$ the latter are given by Kostka numbers, see 
\cite{Sa}. Solution to this problem in a special case would give 
the answer to the famous Foulkes' conjecture, \cite{Fo}.

In Section~\ref{s3} we study some further properties of simple modules over 
factorpowers. We show that there is a contravariant exact involutive
equivalence on the category of $\fp(G,M)$-modules, that all 
simple $\fp(G,M)$-modules are unitarizable and that tensor products
of simple $\fp(G,M)$-modules are completely reducible.
\vspace{2mm}

\noindent
{\bf Acknowledgments.} This work is partially supported by the
Swedish Research Council. The author thanks Ganna Kudryavtseva,
Rowena Paget and Mark Wildon for stimulating discussions.

\section{Irreducible complex representations of the semigroup
$\fp(G,M)$}\label{s2}

\subsection{Green's relations on $\fp(G,M)$}\label{s2.1}

Let $G$ be a finite group acting on a finite set 
$M=\{m_1,m_2,\dots,m_n\}$.  For $A\subset G$ 
let $\overline{A}$ denote the equivalence class of $A$ with respect to 
the equivalence relation $\sim_M$ defining $\fp(G,M)$. For $m\in M$ 
set $A_m=\{\alpha(m)|\alpha\in A\}$. From the definitions for $A,B\subset G$ 
we have that $\overline{A}=\overline{B}$ implies $A_m=B_m$ for all $m$. 
Thus the element $\overline{A}$ of $\fp(G,M)$ can be represented as 
the following array:
\begin{equation}\label{element}
\overline{A}=\left(
\begin{array}{cccc}m_1&m_2&\dots&m_n\\A_{m_1}&A_{m_2}&\dots&A_{m_n}\end{array}
\right).
\end{equation}
Observe that not an arbitrary sequence of subsets of $M$ appears in the 
second row of the above presentation for some elements of $\fp(G,M)$. 
All $A_m$'s are non-empty by definition. Moreover, for any 
$m\in M$ and for any $x\in A_m$ there exists $\sigma\in G$
such that $\sigma(m)=x$ and $\sigma(m')\in A_{m'}$ for all 
$m'\in M$. These two conditions are, obviously, also sufficient. 

Using this notation the multiplication in $\fp(G,M)$ can be written as follows:
{\tiny
\begin{displaymath}
\left(
\begin{array}{ccc}m_1&\dots&m_n\\A_{m_1}&\dots&A_{m_n}\end{array}
\right)\cdot
\left(
\begin{array}{ccc}m_1&\dots&m_n\\B_{m_1}&\dots&B_{m_n}\end{array}
\right)=\left(
\begin{array}{ccc}m_1&\dots&m_n\\ \displaystyle
\bigcup_{x\in B_{m_1}}A_x&\dots&
\displaystyle\bigcup_{x\in B_{m_n}}A_x\end{array}
\right)
\end{displaymath}
}

The anti-involution $\sigma\mapsto\sigma^{-1}$ on $G$ extends in the
natural way to the anti-involution $\pi\mapsto\pi^{\star}$ on $\fp(G,M)$. 
If $\pi=\overline{A}$ for some $A\subset G$, then $\pi^{\star}=
\overline{\{\sigma^{-1}|\sigma\in A\}}$.

The principal result of \cite{M2} asserts that for $\pi,\tau\in\FP$
the condition $\pi\mathcal{L}\tau$ is equivalent to the condition
$\pi=\sigma\tau$ for some $\sigma\in S_n$. Applying $\star$ we have that
$\pi\mathcal{R}\tau$ if and only if $\pi=\tau\sigma$ for some 
$\sigma\in S_n$. We start with extending this result to 
any $\fp(G,M)$.

\begin{proposition}\label{lemgreen}
Let $\pi,\tau\in \fp(G,M)$. Then we have the following:
\begin{enumerate}[(a)]
\item\label{lemgreen.1} $\pi\mathcal{L}\tau$ if and only if there exists
$\sigma\in G$ such that $\pi=\sigma\tau$.
\item\label{lemgreen.2} $\pi\mathcal{R}\tau$ if and only if there exists
$\sigma\in G$ such that $\pi=\tau\sigma$.
\item\label{lemgreen.3} $\pi\mathcal{H}\tau$ if and only if there exists
$\sigma,\sigma'\in G$ such that $\pi=\sigma\tau=\tau\sigma'$.
\item\label{lemgreen.4} $\pi\mathcal{D}\tau$ if and only if there exist
$\sigma,\sigma'\in G$ such that $\pi=\sigma\tau\sigma'$.
\item\label{lemgreen.5} $\mathcal{D}=\mathcal{J}$.
\end{enumerate}
\end{proposition}

\begin{proof}
As $\fp(G,M)$ is finite by definition, the statement \eqref{lemgreen.5}
is well-known, see for example \cite[Proposition~2.1.5]{Ho}. The statements
\eqref{lemgreen.3} and \eqref{lemgreen.4} follow immediately from the
statements \eqref{lemgreen.1} and \eqref{lemgreen.2}. The 
statement \eqref{lemgreen.2} follows from the statement \eqref{lemgreen.1}
applying the anti-involution $\star$. Hence to complete the proof we have to
prove the statement \eqref{lemgreen.1}.
 
Let $A,B\subset G$ be non-empty and assume that $\overline{A}\mathcal{L}
\overline{B}\,\overline{A}$. Let $\tau\in B$ be any element. Since 
$\tau=\overline{\{\tau\}}$ is invertible in $\fp(G,M)$, we have
$\overline{B}\,\overline{A}\mathcal{L}\tau^{-1}\overline{B}\,\overline{A}$.
In particular, $\overline{A}\mathcal{L}
\tau^{-1}\overline{B}\,\overline{A}$ and hence there exists
some $C\subset G$ such that $\overline{A}=\overline{C}\tau^{-1}
\overline{B}\,\overline{A}$. 

Assume that $\overline{A}$ is given by \eqref{element} and 
$\tau^{-1} \overline{B}\,\overline{A}$ is given by a similar formula
with sets $X_{m_i}$'s in the second row. For $B'=\tau^{-1} B$ we observe 
that $B'$ contains the identity element  $e$ of $G$. This means that
$A_{m_i}\subset X_{m_i}$ for all $i$. Assume that $A_{m_i}\subsetneq 
X_{m_i}$ for some $i$, in particular, $|X_{m_i}|>|A_{m_i}|$. 
Write $\overline{C}\tau^{-1} \overline{B}\,\overline{A}$ in the 
form \eqref{element} with sets $Y_{m_i}$'s in the second row, and let
$\sigma\in C$ be arbitrary. Then $|\sigma(X_{m_i})|=|X_{m_i}|$ as
$\sigma$ is an invertible transformation on $M$, and thus
\begin{displaymath}
|Y_{m_i}|\geq |\sigma(X_{m_i})|=|X_{m_i}|>|A_{m_i}|.
\end{displaymath}
This implies that $Y_{m_i}\neq A_{m_i}$ contradicting to 
$\overline{A}=\overline{C}\tau^{-1} \overline{B}\,\overline{A}$.

Thus $A_{m_i}=X_{m_i}$  for all $i$ and hence $\tau^{-1}
\overline{B}\,\overline{A}=\overline{A}$, implying
$\overline{B}\,\overline{A}=\tau\overline{A}$. This completes the proof.
\end{proof}

\subsection{Idempotents and regular $\mathcal{D}$-classes 
in   $\fp(G,M)$}\label{s2.2}

Every subgroup $H$ of $G$ acts on $M$ via restriction. Two subgroups
of $G$ will be called {\em orbit-equivalent} if they have the same
orbits on $M$.

\begin{lemma}\label{lem2.2-1}
Let $H$ be a subgroup of $G$, and $\mathcal{O}(H)$ denote the set of
all subgroups of $G$, which are orbit equivalent to $H$. Then 
$\mathcal{O}(H)$ is partially ordered with respect to inclusion and
contains a unique maximal element, namely, the maximal with respect 
to inclusions set $H'$ in $\overline{H}$.
\end{lemma}

\begin{proof}
If $X$ is an orbit of $H$ on $M$, then every element of $H$ preserves
$X$. Hence, by definition, every element of $H'$ preserves
$X$ as well. On the other hand, the set of all elements from $G$,
which preserve all orbits of $H$ obviously forms a subgroup, say
$\hat{H}$. Finally, if $\sigma\in \hat{H}$, then for any $m\in M$
the element $\sigma(m)\in\{\tau(m)|\tau\in H\}$. Hence
$\hat{H}\subset H'$ and thus $\hat{H}=H'$. The claim follows.
\end{proof}

The maximal element in $\mathcal{O}(H)$ will be called an
{\em orbit-maximal} subgroup of $G$.

\begin{corollary}\label{cor2.2-2}
The idempotents in $\fp(G,M)$ are precisely the classes of 
orbit-maximal subgroups of $G$.
\end{corollary}

\begin{proof}
If $H$ is a subgroup of $G$, then $\overline{H}$ is obviously and
idempotent. From Lemma~\ref{lem2.2-1} we have that $\overline{H}$
is an orbit-maximal subgroup of $G$.

Conversely, let $A\subset G$ be a non-empty subset such that 
it is maximal with respect to inclusions in the class
$\overline{A}$ and assume that $\overline{A}\,\overline{A}=\overline{A}$.
The latter equality means that  $A$ is closed with repsect to 
compositions and hence is a subgroup of $G$ since $G$ is finite. 
The claim follows.
\end{proof}

\begin{remark}\label{rem2.2-3}
{\rm 
If $G=S_n$ with the natural action on $\mathbf{N}$, then the
orbit-maximal subgroups of $G$ are $S_{N_1}\times \dots\times S_{N_k}$,
where $\mathbf{N}=N_1\cup\dots\cup N_k$ is a partition of $\mathbf{N}$
into a disjoint union of nonempty subsets. Thus the idempotents of 
$\FP$ correspond to equivalence relations on $\mathbf{N}$, see
\cite[Theorem~3]{GM2}.
}
\end{remark}

Let $H$ be some orbit-maximal subgroup of $G$ and $\mathcal{D}_{H}$ 
denote the regular $\mathcal{D}$-class of $\fp(G,M)$ containing $H$.
For the study of irreducible representations of $\fp(G,M)$ it is of 
principal importance to understand the structure of 
the {\em trace} $\hat{\mathcal{D}}_{H}=\mathcal{D}_{H}\cup\{0\}$
of $\mathcal{D}_{H}$, which is a semigroup with multiplication
defined as follows:
\begin{displaymath}
\pi\cdot \tau=
\begin{cases}
\pi\tau,& \pi,\tau,\pi\tau\in \mathcal{D}_{H};\\
0,& \text{otherwise}.
\end{cases}
\end{displaymath}
It turns out that the semigroup $\hat{\mathcal{D}}_{H}$ has a very nice
structure.

\begin{proposition}\label{lem2.2-4}
Let $H$ and $H'$ be two orbit-maximal subgroups of $G$. 
\begin{enumerate}[(a)]
\item\label{lem2.2-4.1} We have $H'\in \mathcal{D}_{H}$ if and only if 
$H$ and $H'$ are conjugated in $G$.
\item\label{lem2.2-4.2} The semigroup $\hat{\mathcal{D}}_{H}$ is an
inverse semigroup.
\end{enumerate}
\end{proposition}

\begin{proof}
The ``if'' part of \eqref{lem2.2-4.1} follows directly
from Proposition~\ref{lemgreen}\eqref{lemgreen.4}. To prove the 
``only if'' part assume that $H'\in \mathcal{D}_{H}$. Then, by
Proposition~\ref{lemgreen}\eqref{lemgreen.4}, there exist $\sigma,\tau\in
G$ such that $H'=\sigma H\tau$. As $H'$ has the identity element, we get
$\sigma^{-1}\tau^{-1}\in H$ and hence 
\begin{displaymath}
H'=\sigma H\tau=\sigma (H\sigma^{-1}\tau^{-1})\tau =\sigma H\sigma^{-1},
\end{displaymath}
completing the proof of \eqref{lem2.2-4.1}.

To prove \eqref{lem2.2-4.2} we assume that $H$ and $H'$ are different
orbit-maximal conjugated subgroups of $G$. Then there exists at least
one $H$-orbit which intersects at least two $H'$-orbits and vice versa.
Assume that $H$ is given by \eqref{element} and write $H'H$ in the same
form with subsets $X_{m_i}$'s in the second row. As $H'$ contains the
identity element, we will have $A_{m_i}\subset X_{m_i}$ for all $i$
and $A_{m_i}\neq X_{m_i}$ for at least one $i$. This obviously implies
that $H'H\neq \tau H\sigma$ for any $\tau,\sigma\in G$. Hence
$H'H\not\in \mathcal{D}_{H}$ by Proposition~\ref{lemgreen}\eqref{lemgreen.4},
which yields $H'H=0$ in $\hat{\mathcal{D}}_{H}$. Analogously
$HH'=0$ in $\hat{\mathcal{D}}_{H}$. Hence $\hat{\mathcal{D}}_{H}$ is
a regular semigroup with commuting idempotents, hence inverse
(see for example \cite[Theorem~5.1.1]{Ho}).
\end{proof}

One of the most natural regular $\mathcal{D}$-classes of a semigroup is the
group of units. For $\fp(G,M)$ the latter has the following structure:

\begin{proposition}\label{units}
The mapping $\varphi:\tau\mapsto \overline{\{\tau\}}$, $\tau\in G$, 
is an epimorphism from $G$ to the group of units of $\fp(G,M)$. The
kernel of $\varphi$ coincides with the kernel $K$ of the action of $G$ on $M$. 
\end{proposition}

\begin{proof}
Obviously,  $\varphi$ is a homomorphism from $G$ to the group of units 
of $\fp(G,M)$. The fact that the kernel of this homomorphism coincides 
with $K$ follows from the definition of $\sim_M$. 

Let $e$ denote the identity element of $G$. From the definition of 
$\sim_M$ we have $\overline{e}=K$. Moreover, for any $A\subset G$
and for any $\tau\in A$ we have $\tau K\subset \overline{A}$. If
$A\neq \tau K$ for any $\tau\in A$, then there exist $\tau,\sigma\in G$
and $m\in M$ such that $\tau(m)\neq \sigma(m)$. This means that,
representing $\overline{A}$ in the form \eqref{element}, at least one
of $A_{m_i}$'s in the second row will contain more than one element.
This yields that such $\overline{A}$ cannot be invertible in $\fp(G,M)$,
which implies that $\varphi$ is surjective.
\end{proof}

Proposition~\ref{units} motivates the following statement:

\begin{corollary}\label{faithful}
Let $K$ denote the kernel of the action of $G$ on $M$.
Then the semigroups $\fp(G,M)$ and $\fp(G/K,M)$ are
canonically isomorphic.
\end{corollary}

\begin{proof}
Let $A\subset G$. Then from the proof of  Proposition~\ref{units} we
have that $\overline{A}$ is a union of cosets from $G/K$. It is 
straightforward to verify that mapping $\overline{A}$ to the 
corresponding union of cosets from $G/K$ defines the necessary canonical
isomorphism from the semigroup $\fp(G,M)$ to the semigroup $\fp(G/K,M)$.
\end{proof}

Corollary~\ref{faithful} says that without loss of 
generality in what follows we may assume that the action of $G$ on
$M$ is faithful.

We complete this subsection with a description of maximal subgroups
in $\fp(G,M)$.

\begin{proposition}\label{maxsubgp}
Let $H$ be an orbit-maximal subgroup of $G$ and $\mathrm{N}_G(H)$ its
normalizer in $G$. Then the mapping $\varphi:X\mapsto\overline{X}$
is an isomorphism from $\mathrm{N}_G(H)/H$ to the maximal subgroup
of $\fp(G,M)$ whose identity element is $H$.
\end{proposition}

\begin{proof}
The fact that $\varphi$ is a homomorphism follows from definitions.
The kernel of $\varphi$ consists of all cosets which contain
only those $\tau\in \mathrm{N}_G(H)$ whose action preserves all orbits
of $H$. Since $H$ is orbit maximal, it follows that the kernel of
$\varphi$ coincides with $H$. Thus $\varphi$ is a monomorphism.

Finally, the maximal subgroup of $\fp(G,M)$ with identity $H$ coincides
with the $\mathcal{H}$-class containing $H$. Hence, by 
Proposition~\ref{lemgreen}\eqref{lemgreen.3}, any element from this
maximal subgroup has the form $\overline{\sigma H}=\overline{H\tau}$
for some $\sigma,\tau\in G$. Since $H$ is orbit-maximal, from 
Lemma~\ref{lem2.2-1}  it follows that we thus
must have the equality $\sigma H=H\tau$, that is $H=\sigma^{-1}H\tau$.
This implies $\sigma\tau^{-1}\in H$ and thus $H=\sigma^{-1}H\sigma$, that 
is $\sigma\in \mathrm{N}_G(H)$. This means that $\overline{\sigma H}$
belongs to the image of $\varphi$ and hence $\varphi$ is surjective.
\end{proof}

\subsection{Simple modules over essentially inverse semigroups}\label{s2.3}

A finite semigroup $S$ will be called {\em essentially inverse} provided
that the trace of every regular $\mathcal{D}$-class of $S$ is an
inverse semigroup. For example, if $G$ is a finite group acting on a 
finite set $M$, then the semigroup $\fp(G,M)$ is essentially inverse by
Proposition~\ref{lem2.2-4}. It turns out that for essentially inverse
semigroups the description of simple modules (over complex numbers)
can be substantially simplified in comparison to the classical approach
for all finite semigroups, see \cite{Mu}, \cite{GMS}, 
\cite[Chapter~11]{GM5}. In the following we mostly follow the approach of 
\cite[Chapter~11]{GM5} and \cite{GMS} using the language of bimodules 
and functors.

Let $S$ be an  essentially inverse monoid and $\mathcal{D}_1$,\dots,
$\mathcal{D}_k$ be a complete list of regular $\mathcal{D}$-classes in 
$S$. For every $i=1,\dots,k$ fix some idempotent $e_i\in \mathcal{D}_i$
and let $\mathcal{L}_i$ and $\mathcal{H}_i$ denote the 
$\mathcal{D}$-class and the $\mathcal{H}$-class of $e_i$, respectively.
Then $G_i=\mathcal{H}_i$ is the maximal subgroup of $S$ corresponding
to the idempotent $e_i$. Finally, let $V_i$ denote the formal 
$\mathbb{C}$-span of all elements in $\mathcal{L}_i$. 

\begin{lemma}\label{bimodule}
The assignment
\begin{displaymath}
s\cdot x\cdot g=
\begin{cases}
sxg,& sxg\in \mathcal{L}_i;\\
0,& \text{otherwise}; 
\end{cases}
\end{displaymath}
where $x\in\mathcal{L}_i$, $s\in S$ and $g\in G_i$, defines on 
$V_i$ the structure of an $S-G_i$-bimodule. 
\end{lemma}

\begin{proof}
The left ideal $Se_i$ is obviously stable with respect to the left multiplication with elements from $S$ and right multiplication with 
elements from $G_i$. The same holds for the subset 
$Se_i\setminus \mathcal{L}_i$ of $Se_i$. Now the claim follows from 
the  associativity of the multiplication in $S$.
\end{proof}

Now we are ready to classify all simple $S$-modules.

\begin{theorem}\label{allsimple}
\begin{enumerate}[(a)]
\item\label{allsimple.1} Let $i\in\{1,\dots,k\}$ and $X$ be a simple
$G_i$-module. Then the $S$-module $L(i,X)=V_i\otimes_{G_i}X$ is simple.
\item\label{allsimple.2} Every simple $S$-module has the form
$L(i,X)$ for some $i\in{1,\dots,k}$ and some simple $G_i$-module $X$. 
\item\label{allsimple.3} Modules $L(i,X)$ and $L(j,Y)$ are isomorphic
if and only if $i=j$ and $X\cong Y$.
\end{enumerate}
\end{theorem}

\begin{proof}
Mutatis mutandis \cite[Theorem~11.3.1]{GM5}.
\end{proof}

Denote by $A$ the associative algebra $\oplus_{i=1}^k\mathbb{C}[G_i]$.
Then $V=\oplus_{i=1}^k V_i$ has the natural structure of a
$\mathbb{C}[S]-A$-bimodule. As immediate corollaries from 
Theorem~\ref{allsimple} we have the following two statement:

\begin{corollary}\label{equivalence}
The functor $V\otimes_{A}{}_-$ is an equivalence from the category of
left $A$-modules to the category of all semisimple left 
$\mathbb{C}[S]$-modules. 
\end{corollary}

\begin{proof}
This follows directly from  Theorem~\ref{allsimple}.
\end{proof}

\begin{corollary}\label{jacobson}
Let $J$ denote the Jacobson radical of $\mathbb{C}[S]$, and $n_i$
denote the number of idempotents in $\mathcal{D}_i$, $i=1,\dots,k$. 
Then we have
\begin{displaymath}
\mathbb{C}[S]/J\cong \bigoplus_{i=1}^k \mathrm{Mat}_{n_i}(\mathbb{C}[G_i]),  
\end{displaymath}
where $\mathrm{Mat}_{n_i}(\mathbb{C}[G_i])$ is the algebra of 
$n_i\times n_i$ matrices with coefficients from $\mathbb{C}[G_i]$.
\end{corollary}

\begin{proof}
That the two algebras are Morita equivalent follows from 
Corollary~\ref{equivalence}. It remains to compare the dimensions
of simple modules, which is a straightforward calculation.
\end{proof}

\subsection{Simple modules over factorpowers}\label{s2.4}

Here we combine the results of the two previous subsections to give
a description of simple modules over $\fp(G,M)$. The abstract theory
sounds completely satisfactory for the rough description of these
modules. However, we observe that $G$ is the group of units in $\fp(G,M)$
(note that we already assume that the action of $G$ on $M$ is faithful)
and make the main emphasis on the problem of understanding the structure 
of simple $\fp(G,M)$-modules, when considered as $G$-modules by restriction.

Let $H$ be an orbit-maximal subgroup of $G$. By Theorem~\ref{allsimple},
the simple $\fp(G,M)$-modules corresponding to the idempotent $H$ are
indexed by simple modules over the maximal subgroup of $\fp(G,M)$,
corresponding to $H$. By Proposition~\ref{maxsubgp}, this maximal
subgroup is $\mathrm{N}_G(H)/H$. If $X$ is a simple $\mathrm{N}_G(H)/H$-module,
we denote the corresponding simple $\fp(G,M)$-module by $L(H,X)$.
Denote also by $V_H$ the $\fp(G,M)-\mathrm{N}_G(H)/H$-bimodule
from Subsection~\ref{s2.3}, which corresponds to the idempotent $H$.
Combining the above results we thus obtain:

\begin{theorem}\label{simmodfp}
\begin{enumerate}[(a)]
\item\label{simmodfp.1}
After restriction to $G$ the $G-\mathrm{N}_G(H)/H$-bimodule 
$V_H$  is isomorphic to the bimodule $\mathbb{C}[G/H]$.
\item\label{simmodfp.2}
Let $X$ be a simple $\mathrm{N}_G(H)/H$-module. Then, after restriction to
$G$, the module $L(H,X)$ becomes isomorphic to 
$\mathbb{C}[G/H]\otimes_{\mathrm{N}_G(H)/H}X$.
\end{enumerate}
\end{theorem}

\begin{proof}
The statement \eqref{simmodfp.2} is an immediate corollary from the 
statement \eqref{simmodfp.1} and Theorem~\ref{allsimple}. To prove
\eqref{simmodfp.1} we observe that, by 
Proposition~\ref{lemgreen}\eqref{lemgreen.1}, every element of the
$\mathcal{L}$-class of $H$ has the form $\overline{\tau H}$, $\tau\in G$.
Moreover, since $H$ is orbit-maximal, we have $\overline{\tau H}=H$
if and only if $\tau\in H$. This identifies the $\mathcal{L}$-class of $H$
with $G/H$. It is straightforward to check that the necessary actions
of $G$ and $\mathrm{N}_G(H)/H$ correspond under this identification.
\end{proof}

We note that the $\fp(G,M)$-module structure on $\mathbb{C}[G/H]$
reduces to the $G$-module structure via
Proposition~\ref{lemgreen}\eqref{lemgreen.1}.

\begin{remark}\label{rem2.4-1}
{\rm 
Theorem~\ref{simmodfp} shows another nice property of factorpowers,
namely that on the semigroup level they capture a lot of the intrinsic
group structure of the original group, which in our case is the structure
of certain representations induced from subgroups.  
}
\end{remark}

\begin{remark}\label{rem2.4-2}
{\rm 
The construction in Theorem~\ref{simmodfp} is a special case of the
general fact that if $H$ is a subgroup of $G$, then the set
$G/H$ carries a natural left action of $G$ given by the left multiplication
and the natural right action of $\mathrm{N}_G(H)$ given by the right
multiplication. These two actions of course commute with each other.
}
\end{remark}

\subsection{Simple modules over $\FP$}\label{s2.5}

In the special case of the natural action of the symmetric group $S_n$ 
on $\mathbf{N}$ the ingredients in Theorem~\ref{simmodfp} turn out to
be very classical objects. Therefore we consider this special case in
details.

Let $\rho$ be an equivalence relation on $\mathbf{N}$. Then 
$\rho$ defines a decomposition of $\mathbf{N}$ into a disjoint union
$N_1\cup\dots\cup N_k$ of nonempty subsets. The orbit-maximal subgroup 
of $S_n$ with orbits $N_1,\dots,N_k$ is of course the direct product
$S_{\rho}=S_{N_1}\times\dots\times S_{N_k}$. The normalizer 
$\mathrm{N}_{S_n}(S_{\rho})$ is generated by $S_{\rho}$ and permutations
which arbitrarily permute those blocks of $N_1\cup\dots\cup N_k$
which have the same cardinality. For example, if $|N_1|=\dots=|N_k|=m$,
the group $\mathrm{N}_{S_n}(S_{\rho})$ is isomorphic to the wreath
products $S_m\wr (S_k,\mathbf{N}_k)$. In the general case the group
$\mathrm{N}_{S_n}(S_{\rho})$ is the direct product of such wreath 
products, which correspond to each fixed cardinality of blocks in 
$N_1\cup\dots\cup N_k$. 

The bimodule $\mathbb{C}[S_N/S_{\rho}]$ is by definition nothing else than
the permutation module for $S_n$ associated with the Young subgroup $S_{\rho}$,
see \cite[2.1]{Sa}. If $\lambda$ is the partition of $n$ describing the
cardinalities of blocks in $\rho$, then $\mathbb{C}[S_N/S_{\rho}]$ has
a natural basis consisting of tabloids of shape $\lambda$. The group 
$\mathrm{N}_{S_n}(S_{\rho})/S_{\rho}$ act on $\mathbb{C}[S_N/S_{\rho}]$
by automorphisms permuting rows of the same length in all tabloids.
Assume that for $i=1,\dots,n$  the relation $\rho$ has $k_i$ blocks 
of cardinality $i$. Then
\begin{displaymath}
\mathrm{N}_{S_n}(S_{\rho})/S_{\rho}\cong S_{k_1}\times S_{k_2}\times
\dots\times S_{k_n} 
\end{displaymath}
where we use the convention that $S_0$ is the group with one element. 
Set $\mathbf{k}=(k_1,\dots,k_n)$ and call $\mathbf{l}=(l_1,\dots,l_n)$
a {\em partition} of $\mathbf{k}$ provided that every $l_i$ is a
partition of the corresponding $k_i$ (denoted $\mathbf{l}\vdash\mathbf{k}$).
Then simple $\mathrm{N}_{S_n}(S_{\rho})/S_{\rho}$-modules are just outer tensor
products of Specht modules over the symmetric components and thus 
are indexed by all possible partitions $\mathbf{l}\vdash\mathbf{k}$.
We denote these modules by $S^{\mathbf{l}}$. 

\begin{corollary}\label{corsym}
We have the following decomposition of left $\FP$-modules:
\begin{displaymath}
\mathbb{C}[S_N/S_{\rho}]\cong \oplus_{\mathbf{l}\vdash\mathbf{k}} 
\dim(S^{\mathbf{l}}) L(S_{\rho},S^{\mathbf{l}}).
\end{displaymath}
\end{corollary}

\begin{proof}
The $\FP$-module $\mathbb{C}[S_N/S_{\rho}]$ is obtained via tensor
induction from the regular $\mathrm{N}_{S_n}(S_{\rho})/S_{\rho}$-module.
In the latter each $S^{\mathbf{l}}$ appears with multiplicity 
$\dim(S^{\mathbf{l}})$. Now the claim follows from 
Corollary~\ref{equivalence}.
\end{proof}

Restricting back to $S_n$ we get from Corollary~\ref{corsym} a finer
decomposition of the permutation module into a direct sum of 
submodules. This decomposition is very natural from the $\FP$-point of view. 
It is thus reasonable to ask what is the structure of these
submodules $L(S_{\rho},S^{\mathbf{l}})$, when we consider them 
as $S_n$-modules. 

\begin{problem}\label{bigproblem}
For $\lambda\vdash n$ determine the multiplicity of the Specht module
$S^{\lambda}$ in the $S_n$-module $L(S_{\rho},S^{\mathbf{l}})$.
\end{problem}

For the original permutation module the answer is 
well-known and given by Kostka numbers, see \cite[2.11]{Sa}. Some
special cases of Problem~\ref{bigproblem} are computed in \cite{Ch}.

\begin{remark}\label{remproblem}
{\rm
Problem~\ref{bigproblem} seems to be very difficult in the general case. 
Actually a very special case of it is closely related (and would
give an answer) to the Foulkes' conjecture, see \cite{Fo}. This 
conjecture can be formulated as follows: Let $n=km$, where $k<m$, and 
consider two equivalence relations  $\rho_1$ and $\rho_2$ on $\mathbf{N}$,
the first having $k$ blocks with $m$ elements each, and the second
having $m$ blocks with $k$ elements each. Let $\mathrm{triv}$ 
denote the trivial module
(for any group). Then conjecture says that the multiplicity of
$S^{\lambda}$ in $L(S_{\rho_1},\mathrm{triv})$ does not exceed 
the multiplicity of $S^{\lambda}$ in $L(S_{\rho_2},\mathrm{triv})$.
Although Foulkes' conjecture is proved in some special cases
(see for example \cite{Do,MN,Py} and references therein), the 
general case seems to be wide open.
}
\end{remark}

\subsection{Connection with $\mathcal{F}_n^*$}\label{s2.6}

Let $\mathcal{F}_n^*$ denote the set of all binary relations on
$\mathbf{N}_n$ which have the form $\rho\sigma$, where
$\rho$ is an equivalence relation and $\sigma\in S_n$. Given 
$\rho_1\sigma_1,\rho_2\sigma_2\in \mathcal{F}_n^*$ we define 
\begin{displaymath}
\rho_1\sigma_1\bullet\rho_2\sigma_2=
\overline{\rho_1\sigma_1\rho_2\sigma_1^{-1}}
\sigma_1\sigma_2,
\end{displaymath}
where $\overline{\rho_1\sigma_1\rho_2\sigma_1^{-1}}$ denotes the
minimal equivalence relation, containing the binary relation
$\rho_1\sigma_1\rho_2\sigma_1^{-1}$ (the latter one is the product of
two equivalence relations $\rho_1$ and $\sigma_1\rho_2\sigma_1^{-1}$). 
This makes $\mathcal{F}_n^*$ into an inverse  monoid, which is 
called the {\em maximal factorizable submonoid of the dual symmetric 
inverse monoid} $\mathcal{I}_n^*$, see \cite[Section~3]{FL}. 
Our main result in this subsection is the following:

\begin{theorem}\label{thmdualinv}
Let $J$ denote the Jacobson radical of $\mathbb{C}[\FP]$. Then
$\mathbb{C}[\FP]/J\cong \mathbb{C}[\mathcal{F}_n^*]$.
\end{theorem}

\begin{proof}
Let $\lambda\vdash n$. For $i=1,\dots,n$ denote by $k_{\lambda,i}$ the
number of entries of  $\lambda$, which are equal to $i$. Set
\begin{displaymath}
n_{\lambda}=\frac{n!}{\displaystyle \prod_{i=1}^n k_{\lambda,i}!
\cdot i!^{k_{\lambda,i}}}
\end{displaymath}
and $G_{\lambda}=S_{k_{\lambda,1}}\times S_{k_{\lambda,2}}\times 
\dots\times S_{k_{\lambda,n}}$.

Since $\mathbb{C}[\FP]$ is essentially inverse, using
Corollary~\ref{jacobson} and the description of Green's relations
on $\FP$ (see \cite{M2} or Subsections~\ref{s2.1} and \ref{s2.2}) we have:
\begin{equation}\label{eqeq1}
\mathbb{C}[\FP]/J\cong 
\bigoplus_{\lambda\vdash n} \mathrm{Mat}_{n_{\lambda}}
(\mathbb{C}[G_{\lambda}]).
\end{equation}

Since $\mathcal{F}_n^*$ is an inverse semigroup, the algebra 
$\mathbb{C}[\mathcal{F}_n^*]$ is semi-simple, see for example
\cite[Theorem~11.5.3]{GM5}. In particular, the Jacobson radical of
$\mathbb{C}[\mathcal{F}_n^*]$ is zero. Applying
Corollary~\ref{jacobson} and the description of Green's relations
on $\mathcal{F}_n^*$ (see \cite[Section~3]{FL}) we have:

\begin{equation}\label{eqeq2}
\mathbb{C}[\mathcal{F}_n^*]\cong 
\bigoplus_{\lambda\vdash n} \mathrm{Mat}_{n_{\lambda}}
(\mathbb{C}[G_{\lambda}]).
\end{equation}

The statement now follows by combining \eqref{eqeq1} and \eqref{eqeq2}. 
\end{proof}

As a standard corollary we have:

\begin{corollary}\label{bred}
There exists an algebra monomorphism
$\mathbb{C}[\mathcal{F}_n^*]\hookrightarrow \mathbb{C}[\FP]$.
\end{corollary}

\begin{problem}\label{embproblem}
Construct some explicit  algebra monomorphism
$\mathbb{C}[\mathcal{F}_n^*]\hookrightarrow \mathbb{C}[\FP]$. 
\end{problem}

\begin{remark}\label{comparisonrem}
{\rm  
Theorem~\ref{thmdualinv} basically says that the algebras $\mathbb{C}[\FP]$ 
and $\mathbb{C}[\mathcal{F}_n^*]$ have ``the same'' simple modules. 
This means the following: as the regular  $\mathcal{D}$-classes of $\FP$ and 
the $\mathcal{D}$-classes of $\mathcal{F}_n^*$ are both indexed by partitions 
of $n$, we have an explicit bijection between these $\mathcal{D}$-classes,
which induces a natural bijection between simple modules.
Using the construction of simple modules given in Subsection~\ref{s2.3} 
one sees that the corresponding simple modules over $\mathbb{C}[\FP]$  and 
$\mathbb{C}[\mathcal{F}_n^*]$ are isomorphic when considered
as $S_n$-modules (in particular they have the same dimension).
}
\end{remark}

\section{Duality, unitarizability and tensor products}\label{s3}

Let $G$ be a finite group faithfully acting on a finite set $M$. Set
$A=A(G,M)=\mathbb{C}[\fp(G,M)]$. In this section we will try to get
some information about the category $A\mathrm{-mod}$ of all 
finite-dimensional left $A$-modules.

\subsection{Contravariant duality on the module category}\label{s3.1}

The anti-involution $\star$ on $\fp(G,M)$ extends to an 
anti-involution $\star$ on $A$ in the natural way.
For $V\in A\text{-mod}$ define an $A$-module structure on
$\mathrm{Hom}_{\mathbb{C}}(V,\mathbb{C})$ via $(af)(v):=f(a^{\star}v)$
(note that the fact that $\star$ is an {\em anti}-involution guarantees
that we indeed get a {\em left} $A$-module). This defines a contravariant
functor on $A\text{-mod}$. To distinguish it from the usual duality 
$*$ between $A\text{-mod}$ and $\text{mod-}A$, we will denote the
functor we just defined by $\natural$.

\begin{proposition}\label{duality}
The functor $\natural$ is an exact involutive contravariant equivalence,
which preserves the isomorphism classes of simple $A$-modules.
\end{proposition}

\begin{proof}
That $\natural$ is exact and contravariant follows from the definition.
That $\natural$ is involutive follows from the natural isomorphism 
\begin{displaymath}
 \mathrm{Hom}_{\mathbb{C}}\big(\mathrm{Hom}_{\mathbb{C}}(V,\mathbb{C}),
\mathbb{C}\big)\cong V.
\end{displaymath}
As $\natural\circ\natural$ is isomorphic to the identity functor,
we also get that $\natural$ is an equivalence. Hence we are left to
check that $\natural$ preserves the isomorphism classes of simple $A$-modules.
 
To prove the latter recall from Subsection~\ref{s2.4} that every simple
$A$-module has the form $L(H,X)$, where $H$ is an orbit-maximal subgroup 
of $G$ and $X$ is a simple $\mathrm{N}_G(H)/H$-module. Since $\natural$ is an 
equivalence, the module $L(H,X)^{\natural}$ is simple. Let us determine
this module exactly. First we observe that 
$\pi^{\star}=\pi$ for any idempotent $\pi\in \fp(G,M)$ (as any subgroup 
of $G$ is obviously invariant under taking the inverses of all elements).
From this and the definition of $\natural$ it thus follows that some
idempotent $\pi$ of $\fp(G,M)$ annihilates $L(H,X)^{\natural}$ if and 
only if $\pi$ annihilates $L(H,X)$. This, in particular, implies that 
$L(H,X)^{\natural}$ is associated with the same $\mathcal{D}$-class as
$L(H,X)$. Therefore  $L(H,X)^{\natural}\cong L(H,X')$ for some
simple $\mathrm{N}_G(H)/H$-module $X'$.

To determine the module $X'$ we can compute its character (as
an $\mathrm{N}_G(H)/H$-module) using the definition of $\natural$. First of 
all we observe that the module $X'$ coincides, as a vector space, with the
image of the idempotent $\overline{H}$ (note that the latter is stable
under $\star$). In particular, we get $\dim(X)=\dim(X')$ and can 
restrict our computation to the corresponding subspaces of
$L(H,X)$ and $L(H,X)^{\natural}$. Further, from the definition of $\natural$ 
we see that the action of every  $g\in \mathrm{N}_G(H)/H$ gets substituted by 
the action of $g^{-1}$. However, the  group  $\mathrm{N}_G(H)/H$ is a 
direct product of finite symmetric groups. In every
finite symmetric group the elements $g$ and $g^{-1}$ are conjugated. This
implies that $X$ and $X'$ have the same characters and hence must be
isomorphic. Thus $L(H,X)^{\natural}\cong L(H,X)$, which completes the proof.
\end{proof}

\subsection{Unitarizability of simple modules}\label{s3.2}

The aim of this subsection is to prove that every simple 
$\fp(G,M)$-module is unitarizable with respect to the anti-involution
$\star$ on $\fp(G,M)$ in the ordinary sense. The motivation for this
question was the combination of the two facts: that $\fp(G,M)$ is
essentially inverse (Subsection~\ref{s2.3}), and that simple 
(even all) modules over inverse semigroups are unitarizable
(\cite[11.5]{GM5}). As $\fp(G,M)$ is not inverse in general, the algebra
$A$ is not semi-simple and hence all $A$-modules cannot be unitarizable
in the ordinary sense. However, one could expect that at least simple
modules are. We will prove this in the present subsection.

Let $H$ be an orbit-maximal subgroup of $G$ and $X$ be a simple
$\mathrm{N}_G(H)/H$-module. Let further $H_1=H$, $H_2$,\dots, $H_k$ be a
list of all pairwise differnet subgroups of $G$, conjugated to $H$.
For $i=1,\dots,k$ set $\varepsilon_i=\overline{H_i}$. Then from 
Proposition~\ref{lemgreen}\eqref{lemgreen.4} it follows that 
$\varepsilon_1$,\dots, $\varepsilon_k$ is a list of all pairwise 
differnet idempotents in the regular $\mathcal{D}$-class of $\fp(G,M)$,
containing the idempotent $\overline{H}$. Finally, for $i=2,\dots,k$
fix some $g_i\in G$ such that $H_i=g_i^{-1} H g_i$, and set $g_1=e$.

\begin{lemma}\label{lem3.2-1}
For $i=1,\dots,k$ let $V_i$ denote the image of $\varepsilon_i$ on
$L(H,X)$. Then $L(H,X)=\oplus_{i=1}^k V_i$.
\end{lemma}

\begin{proof}
As $\hat{\mathcal{D}}_H$ is an inverse semigroup by 
Proposition~\ref{lem2.2-4}\eqref{lem2.2-4.2}, the product of any two
different idempotents in it is zero.  If $v_i\in V_i$ and 
$\lambda_i\in\mathbb{C}$ are such that $\sum_i \lambda_i v_i=0$, 
applying $\varepsilon_j$ and using the above we get $\lambda_j v_j=0$
and thus for every $j$ we either have $\lambda_j=0$ or $v_j=0$. This 
means that the sum of all the $V_i$'s  is direct.

On the other hand, let $0\neq v\in L(H,X)$, $v_i= \varepsilon_i v$ and
$v'=\sum_i v_i$. Then $\varepsilon_i(v-v')=0$ for all $i$. Set
\begin{displaymath}
W=\{w\in L(H,X)\vert\varepsilon_i w=0\text{ for all }i\}.
\end{displaymath}
Let  $w\in W$, $\pi\in \fp(G,M)$ and $i\in\{1,\dots,k\}$. If 
$\varepsilon_i \pi\not\in \mathcal{D}_H$, we have 
$\varepsilon_i \pi w=0$ by definition. If $\varepsilon_i \pi\in 
\mathcal{D}_H$, then the left class of $\varepsilon_i \pi$ contains
some $\varepsilon_j$ and hence $\varepsilon_i \pi=
\tau\varepsilon_j$ for some $\tau\in G$ by 
Proposition~\ref{lemgreen}\eqref{lemgreen.1}. This again implies that
$\varepsilon_i \pi w=0$. Thus $W$ is a submodule of $L(H,X)$, different
from $L(H,X)$ (as $\varepsilon_1\neq 0$). Since $L(H,X)$ is simple, we
get $W=0$ and hence $v-v'=0$. This proves that $L(H,X)=\sum_{i=1}^k V_i$
and completes the proof of the lemma.
\end{proof}

Let $(\cdot,\cdot)_1$ be any Hermitian scalar product on $V_1$. For
$v,w\in V_1$ set 
\begin{displaymath}
(v,w)=\sum_{g\in \mathrm{N}_G(H)/H}(g(v),g(w))_1 
\end{displaymath}
and note that this new Hermitian scalar product on $V_1$ is invariant
with respect to the action of $ \mathrm{N}_G(H)/H$.
For $i\neq j$ and arbitrary $v\in V_i$, $w\in V_j$ we set $(v,w)=0$.
Finally, for $v\in V_i$ we note that $\varepsilon_i=g_i^{-1}\varepsilon_1
g_i$ by our choice of $g_i$ and hence  $g_i(v)\in V_1$. Thus for 
$v,w\in V_i$ we may set $(v,w)=(g_i(v),g_i(w))$ and extend
$(\cdot,\cdot)$ to the whole $L(H,X)$ by skew-bilinearity. Now we are in
position to show that the $\fp(G,M)$-module $L(H,X)$ is {\em unitarizable}
in the following sense:

\begin{proposition}\label{unitary}
$(\cdot,\cdot)$ is a Hermitian scalar product on $L(H,X)$ and
\begin{equation}\label{equnitary}
(\sigma(v),w)=(v,\sigma^{\star}(w))
\end{equation}
for all $\sigma\in \fp(G,M)$ and $v,w\in L(H,X)$.
\end{proposition}

\begin{proof}
The product $(\cdot,\cdot)$ is bilinear and skew-symmetric by construction
and its restriction to every $V_i$ is positive definite as 
$(\cdot,\cdot)_1$ is positive definite. Hence the fact that 
$(\cdot,\cdot)$ is positive definite follows from Lemma~\ref{lem3.2-1}.
This shows that $(\cdot,\cdot)$ is a Hermitian scalar product on $L(H,X)$.

Let now $\sigma\in \fp(G,M)$, $v\in V_i$ and $w\in V_j$. If
$\varepsilon_j\sigma\varepsilon_i\not\in \mathcal{D}_H$, then
$(\sigma(v),w)=(v,\sigma^{\star}(w))=0$ follows immediately from the definitions and thus \eqref{equnitary} holds. Hence we may assume 
$\sigma'=\varepsilon_j\sigma\varepsilon_i\in \mathcal{D}_H$ and
even more that $\sigma=\sigma'$. Then $\sigma(v)\in V_j$ and hence we
have $(\sigma(v),w)=(g_j\sigma(v),g_j(w))$. From 
$\sigma=\varepsilon_j\sigma\varepsilon_i\in \mathcal{D}_H$ it follows that
$g_j\sigma g_i^{-1}\in \mathrm{N}_G(H)/H$. Using the invariance of
$(\cdot,\cdot)$ with respect to the action of $\mathrm{N}_G(H)/H$ on
$V_1$ we  get $(g_j\sigma(v),g_j(w))=(g_i(v),g_i\sigma^{\star}(w))$,
which implies the necessary equality \eqref{equnitary} and completes the
proof.
\end{proof}

\begin{remark}
{\rm The existence and uniqueness (in the correct sense) of some 
(not necessarily positive definite) scalar product on $L(H,X)$ having 
the property  \eqref{equnitary} follows easily from 
Subsection~\ref{s3.1} and \cite[Theorem~1]{MT}. However, the fact 
that this scalar product is positive definite (established in 
Proposition~\ref{unitary}) will be crucial for the next subsection.
}
\end{remark}

\subsection{Complete reducibility of tensor products}\label{s3.3}

If $S$ is a semigroup and $X$ and $Y$ are two $S$-modules, then the
tensor product $X\otimes_{\mathbb{C}}Y$ has a natural structure of an
$S$-module via the diagonal action $s(x\otimes y)=s(x)\otimes s(y)$.
Tensor products of certain modules over various generalizations of the
symmetric group appear in the literature, especially in connection to
Schur-Weyl dualities (see \cite{So,KM}).  The main result of this 
subsection is the following statement:

\begin{theorem}\label{tensor}
Let $H$ and $H'$ be two orbit-maximal subgroups of $G$ and
$X$ and $X'$ be two simple modules over $\mathrm{N}_G(H)/H$ and
$\mathrm{N}_G(H')/H'$ respectively. Then the $\fp(G,M)$-module
$L(H,X)\otimes_{\mathbb{C}}L(H',X')$ is completely reducible.
\end{theorem}

\begin{proof}
By Proposition~\ref{unitary}, both $L(H,X)$ and $L(H',X')$ are unitarizable.
Let $(\cdot,\cdot)$ and $(\cdot,\cdot)'$ denote the corresponding Hermitian
scalar products. Let further $v_1,\dots,v_k$ and $w_1,\dots,w_m$ be
some orthonormal bases in $L(H,X)$ and $L(H',X')$ with respect to 
$(\cdot,\cdot)$ and $(\cdot,\cdot)'$ respectively. Let 
$\langle\cdot,\cdot\rangle$ denote the Hermitian
scalar product on the vector space $L(H,X)\otimes_{\mathbb{C}}L(H',X')$ 
for which  $\{v_i\otimes w_j\}$ is an orthonormal basis. 

For $\sigma\in \fp(G,M)$ set $\sigma_{i,j}=(\sigma(v_i),v_j)$ and
$\sigma'_{a,b}=(\sigma(w_a),w_b)'$. Then unitarizability of 
$L(H,X)$ and $L(H',X')$ means that 
$\sigma_{i,j}=\overline{(\sigma^{\star})_{j,i}}$ and
$\sigma'_{a,b}=\overline{(\sigma^{\star})'_{b,a}}$ for all
$\sigma,\sigma',i,j,a,b$, where $\overline{\cdot}$ denotes the 
complex conjugation. Using this we have
\begin{displaymath}
\langle\sigma(v_i\otimes w_a),v_j\otimes w_b\rangle =
\sigma_{i,j}\sigma'_{a,b}
\end{displaymath}
for all appropriate $i,j,a,b$ and, moreover, 
\begin{displaymath}
\sigma_{i,j}\sigma'_{a,b}=\overline{(\sigma^{\star})_{j,i}}
\cdot \overline{(\sigma^{\star})'_{b,a}}=
\overline{(\sigma^{\star})_{j,i}(\sigma^{\star})'_{b,a}},
\end{displaymath}
which implies that the equality \eqref{equnitary} holds for
the form $\langle\cdot,\cdot\rangle$ on 
$L(H,X)\otimes_{\mathbb{C}}L(H',X')$. In particular, the module
$L(H,X)\otimes_{\mathbb{C}}L(H',X')$ is unitarizable. 

Now if $W\subset L(H,X)\otimes_{\mathbb{C}}L(H',X')$ is a subspace invariant
under the action of $\fp(G,M)$, using \eqref{equnitary} one shows that
the orthogonal complement $W^{\perp}$ is also invariant (and is really a 
complement as our scalar product is positive definite). This implies
that $L(H,X)\otimes_{\mathbb{C}}L(H',X')$ is completely reducible.
\end{proof}

The following natural question might be very difficult since, as far as 
the author knows, the complete  answer  to it is not yet known even for 
the symmetric group. However, a satisfactory ``semigroup'' answer could 
be a reduction of this problem to the analogous problem for the 
symmetric group.

\begin{problem}\label{prtensor}
Determine the multiplicity of each $L(H,X)$ in the decomposition of  
$L(H',X')\otimes_{\mathbb{C}}L(H'',X'')$.
\end{problem}

\begin{remark}\label{remend1}
{\rm
Tensoring with $L(H,X)$ defines an exact endofunctor on 
$\fp(G,M)\text{-mod}$. In would be interesting to understand this functor
and its adjoints.  
}
\end{remark}

\begin{remark}\label{remend2}
{\rm
Proposition~\ref{unitary} and Theorem~\ref{tensor} generalize (with the same
proof) to arbitrary essentially inverse semigroup with involution, which
induces the usual involution $a\mapsto a^{-1}$ on all
inverse $\mathcal{D}$-classes.
}
\end{remark}

\vspace{1cm}


\noindent
Department of Mathematics, University of Glasgow,
University Gardens, Glasgow, G12 8QW, UK,
e-mail: {\tt mazor\symbol{64}maths.gla.ac.uk}
\vspace{0.5cm}

\noindent
and Department of Mathematics, Uppsala University, SE 471 06,
Uppsala, SWEDEN, e-mail: {\tt mazor\symbol{64}math.uu.se}

\end{document}